\documentclass[10pt]{article}

\makeatletter
\@addtoreset{equation}{section}
\makeatother

\input epsf

\usepackage{amssymb}
\usepackage{amsmath}
\usepackage{epsfig}

\newtheorem{theorem}{Theorem}[section]

\newtheorem{corollary}[theorem]{Corollary}

\def\tfrac#1#2{{{\lower.6ex
\hbox{$\scriptstyle#1$}}\over
{\raise.7ex
\hbox{$\scriptstyle#2$}}}}

\def\sign{{\rm sign}}

\def\sign{{\rm sign}}
\def\erfc{{\rm erfc}}

\def\CC{{\mathbb C}}

\def\dsp#1{\displaystyle#1}
\def\sign{{\rm sign}}
\def\erfc{{\rm erfc}}

\def\CC{{\mathbb C}}

\def\dsp#1{\displaystyle#1}

\begin{document}
\title{Asymptotic analysis of the normal inverse \\ Gaussian cumulative distribution}
\author{Nico M. Temme\footnotemark[1]}

\date{\today}

\maketitle

\renewcommand{\thefootnote}{\fnsymbol{footnote}}

\footnotetext[1]{Valkenierstraat 25, 1825BD Alkmaar, The Netherlands.  Email: nic@temme.net}

\begin{abstract}

Using a recently derived integral in terms of elementary functions, we derive new asymptotic expansions of the normal inverse Gaussian cumulative distribution function. One of the asymptotic representations is in terms of the normal Gaussian distribution or complementary error function.
\end{abstract}



\section{Introduction}\label{sec:intro}
The normal inverse Gaussian distribution is a four-parameter distribution ($\alpha, \beta,\mu,\delta$) with argument $x$, which has been introduced by Barndorﬀ-Nielsen \cite{Barndorff:1977:EDD}, \cite{Barndorff:1978:HDD}, \cite{Barndorff:1997:NIG}.

 In a recent preprint \cite{Guillermo:2025:NIG}  the commonly used representation of the cumulative distribution function is given by
  \begin{equation}\label{eq:intro01}
F(x;\alpha,\beta,\mu,\delta)=\frac{\alpha\delta e^{\delta\gamma}}{\pi}\int_{-\infty}^x\frac{K_1\left(\alpha\sqrt{\delta^2+(t-\mu)^2}\right)}{\sqrt{\delta^2+(t-\mu)^2}}e^{\beta(t-\mu)}\,dt,
\end{equation} 
where $\gamma=\sqrt{\alpha^2-\beta^2}$ and $K_1(z)$ denotes the modified Bessel function. The cited paper gives new convergent series and derived asymptotic expansions, with the aim  of developing a software package to compute the cumulative distribution function based on the normal inverse Gaussian distribution.

In this paper we give more asymptotic expansions  of $F(x;\alpha,\beta,\mu,\delta)$, after writing this function in a standard form to apply Laplace's method, with a modification to handle the case  that a pole is near a saddle point. This gives an expansion in which the complementary error function controls this phenomenon. 

The starting point of our approach is an integral representation in terms of elementary functions. In fact this is given for the function that we like to call the complementary function, defined by 
 \begin{equation}\label{eq:intro02}
G(x;\alpha,\beta,\mu,\delta)=1-F(x;\alpha,\beta,\mu,\delta).
\end{equation} 
Especially for numerical computations, it is important to have stable representation for both functions, and to compute the smaller of the two first. As we will see later in this paper, the transition value $x_0$ with respect to $x$ is given by 
\begin{equation}\label{eq:intro03}
x_0=\mu+\frac{\beta\delta}{\sqrt{\alpha^2-\beta^2}},
\end{equation} 
and, roughly speaking, because this follows from the asymptotic approximation, we have $0\le F(x;\alpha,\beta,\mu,\delta)\lessapprox\frac12$ when $x\lessapprox x_0$. 

The integral representation we use is for the $G$-function is recently derived  in \cite[Eqn.~(2.8)]{Guillermo:2025:NIG}):
\begin{equation}\label{eq:intro04}
G(x;\alpha,\beta,\mu,\delta)=\frac{e^{\delta\gamma}}{\pi}
\int_0^\infty \frac{r e^{-\xi\left(\sqrt{r^2+\alpha^2}-\beta\right)}}{\sqrt{r^2+\alpha^2}\left(\sqrt{r^2+\alpha^2}-\beta\right)}\sin(\delta r)\,dr,
\end{equation} 
where
\begin{equation}\label{eq:intro05}
\xi=x-\mu\ge0, \quad \delta>0, \quad \alpha>0, \quad -\alpha<\beta<\alpha,\quad \gamma=\sqrt{\alpha^2-\beta^2}.
\end{equation} 
As observed in \cite{Guillermo:2025:NIG}, for $\xi<0$ we can use the relation
\begin{equation}\label{eq:intro06}
\begin{array}{r@{\,}c@{\,}l}
F(x;\alpha,\beta,\mu,\delta)&=&1-F(-x;\alpha,-\beta,-\mu,\delta)=G(-x;\alpha,-\beta,-\mu,\delta)\\[8pt]
&=&\dsp{\frac{e^{\delta\gamma}}{\pi}
\int_0^\infty \frac{r e^{\xi\left(\sqrt{r^2+\alpha^2}+\beta\right)}}{\sqrt{r^2+\alpha^2}\left(\sqrt{r^2+\alpha^2}+\beta\right)}\sin(\delta r)\,dr.}
\end{array}
\end{equation}

In our new asymptotic results, the key term in the approximations is the complementary error function, defined by 
\begin{equation}\label{eq:intro07}
\erfc(z)=\frac{2}{\sqrt{\pi}}\int_{z}^\infty e^{-t^2}\,dt,\quad z\in\CC,
\end{equation} 
and a comparable representation is not available in \cite{Guillermo:2025:NIG}. As we have explained previously, for example in \cite{Gil:2014:TAN}, \cite{Gil:2019:TCI}, \cite{Gil:2022:NAR} and also in \cite[Chapter~21 and Part~7]{Temme:2015:AMI}, such a representation can yield a powerful asymptotic approximation also with respect to so-called uniformity parameters, and moreover it provides an excellent starting point for inverting cumulative distribution functions with respect to one of the parameters. This topic is planned for future research.

In Section~\ref{sec:trans}  we use several transformations of the integral in \eqref{eq:intro04} and obtain in Section~\ref{sec:prep}  a representation suitable for asymptotic analysis. In  Section~\ref{sec:asy} we give the asymptotic expansions, together with a figure and a table to explain the role of the transition value $x_0$. 
In an Appendix we give a short Maple code for the evaluation of the coefficients used in certain asymptotic expansions.

\section{Transformations of the integral}\label{sec:trans}

The integrand in \eqref{eq:intro04} is an even function of $r$ and 
we substitute  $r=\alpha\sinh(t)$. 
 Then we obtain
\begin{equation}\label{eq:trans01}
G(x;\alpha,\beta,\mu,\delta)=\frac{e^{\delta\gamma+\beta\xi}}{2\pi}\int_{-\infty}^\infty e^{-\xi\alpha\cosh(t)}
\frac{\sinh(t)\sin\left(\alpha\delta\sinh(t)\right)}{\cosh(t)-\cos(\tau)}\,dt,
\end{equation} 
where $\tau$ follows from 
\begin{equation}\label{eq:trans02}
\frac{\beta}{\alpha}=\cos(\tau), \quad 0<\tau<\pi.
\end{equation} 
We use $\sin(z)=\frac{1}{i}\left(e^z-\cos(z)\right)$, with $z=\alpha\delta\sinh(t)$, observe that the cosine term will give an odd integrand, and write \eqref{eq:trans01} as
\begin{equation}\label{eq:trans03}
G(x;\alpha,\beta,\mu,\delta)=\frac{e^{\delta\gamma+\beta\xi}}{2\pi i}\int_{-\infty}^\infty e^{-\alpha\omega\phi(t)}
\frac{\sinh(t)}{\cosh(t)-\cos(\tau)}\,dt,
\end{equation} 
where
\begin{equation}\label{eq:trans04}
\phi(t)=\frac{\xi}{\omega}\cosh(t)-i\frac{\delta}{\omega}\sinh(t),\quad \omega=\sqrt{\xi^2+\delta^2}.
\end{equation} 
We introduce $\nu\in(0,\frac12\pi)$ by writing
\begin{equation}\label{eq:trans05}
 \nu=\arctan\frac{\delta}{\xi}\quad \Longrightarrow \quad \xi=\omega\cos(\nu)\quad {\rm and} \quad \delta=\omega\sin(\nu).
\end{equation} 
It follows that the function $\phi(t)$ can be written as
\begin{equation}\label{eq:trans06}
\phi(t)=\cosh(t)\cos(\nu)-i\sinh(t)\sin(\nu)=\cosh(t-i\nu).
\end{equation} 
This gives for the integral representation in \eqref{eq:trans01}
\begin{equation}\label{eq:trans07}
G(x;\alpha,\beta,\mu,\delta)= \frac{e^{\delta\gamma+\beta\xi}}{2\pi i}\int_{-\infty}^\infty e^{-\alpha\omega\cosh(t-i\nu)}
\frac{\sinh(t)}{\cosh(t)-\cos(\tau)}\,dt.
\end{equation} 
This integral converges when $\cos(\nu)\ge0$, but for the representation of $G(x;\alpha,\beta,\mu,\delta)$ in  \eqref{eq:intro04} we assumed that $\delta=\omega\sin(\nu)>0$, and therefore  we assume $\nu\in \left(0,\frac12\pi\right)$.

We want to integrate the integral in \eqref{eq:trans07} along the horizontal path with $\Im t=\nu$, and we need information about the poles of the integrand. Since $-1<\cos(\tau)<1$ the poles closest to the origin are $t_\pm =\pm i\tau$, we can shift the path of integration in \eqref{eq:trans07}  to the path $\Im t=\nu$. When $\nu >\tau$ we cross the pole at $ i\tau$, and calculate the residue. 

 After shifting the path, we integrate along the horizontal line $\Im t=\nu$ using the substitution $t=s+i\nu$, and when  $\nu >\tau$ we evaluate the relation
\begin{equation}\label{eq:trans08}
-\alpha\omega\cosh(i\tau-i\nu)=-\alpha\omega\cos(\tau-\nu)=-\beta\xi-\delta\gamma.
\end{equation}
We obtain the representations
\begin{equation}\label{eq:trans09}
\begin{array}{r@{\,}c@{\,}l}
G(x;\alpha,\beta,\mu,\delta)&=&\dsp{\frac{e^{\delta\gamma+\beta\xi}}{2\pi i}\int_{-\infty}^\infty e^{-\alpha\omega\cosh(s)}
\frac{\sinh(s+i\nu)}{\cosh(s+i\nu)-\cos(\tau)}\,ds,\quad \tau > \nu,}\\[8pt]
G(x;\alpha,\beta,\mu,\delta)&=&\dsp{1-\frac{e^{\delta\gamma+\beta\xi}}{2\pi i}\int_{-\infty}^\infty e^{-\alpha\omega\cosh(s)}
\frac{\sinh(s+i\nu)}{\cos(\tau)-\cosh(s+i\nu)}\,ds,\quad \tau < \nu.}
\end{array}
\end{equation} 
Moreover, the case $\tau<\nu$ gives
\begin{equation}\label{eq:trans10}
F(x;\alpha,\beta,\mu,\delta)=\frac{e^{\delta\gamma+\beta\xi}}{2\pi i}\int_{-\infty}^\infty e^{-\alpha\omega\cosh(s)}
\frac{\sinh(s+i\nu)}{\cos(\tau)-\cosh(s+i\nu)}\,ds,\quad \tau < \nu.
\end{equation}

For convergence of the integral in  \eqref{eq:intro04} we assumed $\xi\ge0$, or $x\ge\mu$, but in the above three integrals this is no longer needed, and we can assume $\xi\le0$, or $x\le\mu$.

These integrals have a saddle point at the origin, and the real axis is the path of steepest descent. There are poles at $s$-values corresponding with the poles $t_\pm=\pm i\tau$, and by the transformation $t=s+i\nu$ the poles in the $s$-variable are given by
\begin{equation}\label{eq:trans11}
 s_\pm=-i(\nu\mp\tau),\quad 
\nu\in(0,\pi),\quad \tau\in(0,\pi). 
\end{equation} 

We see that the pole $s_+$ coincides with the saddle point at the origin when $\xi$ takes the value $\xi_0$ that follows from
\begin{equation}\label{eq:trans12}
\tau=\nu\quad \Longrightarrow\quad \arctan\frac{\delta}{\xi_0}=\arccos\frac{\beta}{\alpha} \quad \Longrightarrow \quad \xi_0=\frac{\delta\beta}{\sqrt{\alpha^2-\beta^2}}=\frac{\delta\beta}{\gamma}.
\end{equation} 
This gives for $x=\mu+\xi$ the transition value $x_0=\xi_0+\mu$ announced in \eqref{eq:intro03}. 
The transition value $x_0 = \xi_0 + \mu$ coincides with the mean of the normal Gaussian distribution.

We have the following cases:
\begin{equation}\label{eq:trans13}
\begin{array}{llll}
&\tau>\nu \quad \Longrightarrow \quad &x> x_0, \quad \Im s_+>0, \quad &\Im s_-<0,\\
&\tau=\nu \quad \Longrightarrow \quad &x= x_0, \quad  s_+=0, \quad &\Im s_-<0,\\
&\tau<\nu \quad \Longrightarrow \quad &x< x_0, \quad \Im s_+<0, \quad &\Im s_-<0.
\end{array}
\end{equation}
When one of these poles is close to the origin, we need special asymptotic methods to deal with it.

\section{Further preparations for the asymptotic analysis}\label{sec:prep}

We continue with \eqref{eq:trans10}, assuming that $\nu>\tau$; this condition will be relaxed after deriving the asymptotic expansions in Section~\ref{sec:asy}.  We write  the integrand in real and imaginary parts, using
\begin{equation}\label{eq:prep01}
\begin{array}{ll}
&\sinh(s+i\nu)=\sinh(s)\cos(\nu)+i\cosh(s)\sin(\nu),\\
&\cosh(s+i\nu)=\cosh(s)\cos(\nu)+i\sinh(s)\sin(\nu).
\end{array}
\end{equation} 
We obtain
\begin{equation}\label{eq:prep02}
\begin{array}{lll}
\dsp{\frac{\sinh(s+i\nu)}{\cos(\tau)-\cosh(s+i\nu)}}&=&\dsp{\frac{\sinh(s)\left(\cos(\tau)\cos(\nu)-\cosh(s)\right)}{\left(\cosh(s)-\cos(\tau)\cos(\nu)\right)^2-\sin^2(\tau)\sin^2(\nu)}}\ +\\[8pt]
&&\dsp i{\frac{\sin(\nu)\left(\cos(\tau)\cosh(s)-\cos(\nu)\right)}
{\left(\cosh(s)-\cos(\tau)\cos(\nu)\right)^2-\sin^2(\tau)\sin^2(\nu)}}. 
\end{array}
\end{equation} 
We see that the real part is odd with respect to $s$ and the imaginary part is even. So we can write:
\begin{equation}\label{eq:prep03}
F(x;\alpha,\beta,\mu,\delta)=\frac{e^{\delta\gamma+\beta\xi}}{2\pi}\int_{-\infty}^\infty 
\frac{e^{-\alpha\omega\cosh(s)}\sin(\nu)\left(\cos(\tau)\cosh(s)-\cos(\nu)\right)}
{\left(\cosh(s)-\cos(\tau)\cos(\nu)\right)^2-\sin^2(\tau)\sin^2(\nu)}\,ds.
\end{equation} 
Next we use  $\cosh(s)=1+2\sinh^2(s/2)$ and substitute
\begin{equation}\label{eq:prep04}
\sigma=\sinh(s/2), \quad  \frac{ds}{d\sigma}=\frac{2}{\sqrt{1+\sigma^2}}.
\end{equation} 
This gives
\begin{equation}\label{eq:prep05}
F(x;\alpha,\beta,\mu,\delta)=\frac{e^{\delta\gamma+\beta\xi-\alpha\omega}}{4\pi}\int_{-\infty}^\infty 
e^{-2\alpha\omega\sigma^2}f(\sigma)\,d\sigma,\quad \nu>\tau,
\end{equation} 
where
\begin{equation}\label{eq:prep06}
\begin{array}{r@{\,}c@{\,}l}
f(\sigma)&=&\dsp{\frac{\sin(\nu)}{\sqrt{1+\sigma^2}}\ \frac{\cos(\tau)(1+2\sigma^2)-\cos(\nu)}
{\left(\sigma^2+\frac12-\frac12\cos(\tau)\cos(\nu)\right)^2-\frac14\sin^2(\tau)\sin^2(\nu)} }\\[8pt]
&=&\dsp{\frac{\sin(\nu)}{\sqrt{1+\sigma^2}}\ \frac{\cos(\tau)(1+2\sigma^2)-\cos(\nu)}{(\sigma^2-\sigma_+^2)(\sigma^2-\sigma_-^2)}.}
\end{array}
\end{equation}

The poles $\sigma_\pm$ follow from the poles $s_\pm$ via the relation $\sigma=\sinh(s/2)$. We have (see \eqref{eq:trans11})
\begin{equation}\label{eq:prep07}
\sigma_\pm=\sinh\left(\tfrac12s_\pm\right)=i\sin\left(\tfrac12s_\pm\right)=-i\sin\left(\tfrac12(\nu\mp\tau)\right).
\end{equation} 
By splitting into fractions we obtain
\begin{equation}\label{eq:prep08}
\begin{array}{r@{\,}c@{\,}l}
f(\sigma)&=&\dsp{\frac{\sin(\nu)}{\sqrt{1+\sigma^2}\,(\sigma_-^2-\sigma_+^2)}
\left(\frac{\cos(\nu)-\cos(\tau)(1+2\sigma_+^2)}{\sigma^2-\sigma_+^2}-\frac{\cos(\nu)-\cos(\tau)(1+2\sigma_-^2)}{\sigma^2-\sigma_-^2}\right)}\\[8pt]
&=&\dsp{\frac{1}{ \sqrt{1+\sigma^2}}
\left(\frac{\sin(\nu-\tau)}{\sigma^2-\sigma_+^2}+\frac{\sin(\nu+\tau)}{\sigma^2-\sigma_-^2}\right)}.
\end{array}
\end{equation}

The argument of the exponential function in front of the integral in \eqref{eq:prep05} can be written  as
\begin{equation}\label{eq:prep09}
\delta\gamma+\beta\xi-\alpha\omega=-2\alpha\omega\sin^2\left(\tfrac12(\nu-\tau)\right)=2\alpha\omega\sigma_+^2.
\end{equation}

After these steps, we summarize the results obtained so far in the following theorem.

\begin{theorem}\label{thm:thm01}
We can write \eqref{eq:prep05} in the form
\begin{equation}\label{eq:prep10}
\begin{array}{lll}
&&\dsp{F(x;\alpha,\beta,\mu,\delta)=F^+(x;\alpha,\beta,\mu,\delta)+F^-(x;\alpha,\beta,\mu,\delta),}
\\[8pt]
&&F^\pm(x;\alpha,\beta,\mu,\delta)=\dsp{\frac{e^{z\sigma_+^2}}{4\pi}
\sin(\nu\mp\tau)U(\sigma_\pm,z),}\\[8pt]
&&\dsp{
U(\sigma_\pm,z)=\int_{-\infty}^\infty e^{-z\sigma^2}\frac{d\sigma}{(\sigma^2-\sigma_\pm^2)\sqrt{1+\sigma^2}},}
\end{array}
\end{equation}
where, including the other  notations used so far:
\begin{equation}\label{eq:prep11}
\begin{array}{lll}
&& \dsp{z=2\alpha\omega,\quad \xi=x-\mu, \quad x_0=\mu+\xi_0,\quad \xi_0=\frac{\beta\delta}{\gamma},}\\[8pt]
&&\sigma_+=\dsp{-i\sin\left(\tfrac12(\nu-\tau)\right),\quad \sigma_-=-i\sin\left(\tfrac12(\nu+\tau)\right),}\\[8pt]
&&\dsp{\xi=x-\mu=\omega\cos(\nu),\quad \delta=\omega\sin(\nu),}\\[8pt]
&&\dsp{\nu=\arctan\frac{\delta}{\xi}, \quad \nu\in \left(0,\tfrac12\pi\right),}\\[8pt]
&&\dsp{\omega=\sqrt{\xi^2+\delta^2},\quad\cos(\tau)= \frac{\beta}{\alpha},\quad \tau\in(0,\pi),}\\[8pt]
&&\dsp{\gamma=\sqrt{\alpha^2-\beta^2}=\alpha\sin(\tau).}\\[8pt]
\end{array}
\end{equation} 
\end{theorem}

\section{Asymptotic expansions}\label{sec:asy}
First, we give the  asymptotic expansion of $U(\sigma_\pm,z)$ as defined in \eqref{eq:prep10}, with $z$ as a positive large parameter and $i\sigma_\pm\in(0,1)$. We first assume that  $i\sigma_\pm$ is not small, say $\frac12\le i\sigma_\pm<1$. In that case, an asymptotic expansion can be obtained by using Laplace's method, see \cite[Chapter~3]{Temme:2015:AMI}.

We expand 
\begin{equation}\label{eq:asy01}
\frac{1}{(\sigma^2-\rho^2)\sqrt{1+\sigma^2}}=\sum_{k=0}^\infty u_k(\rho) \sigma^{2k},\quad \rho=\sigma_\pm,
\end{equation} 
and obtain the asymptotic expansion
\begin{equation}\label{eq:asy02}
U(\rho,z)\sim\sqrt{\frac{\pi}{z}}\sum_{k=0}^\infty u_k(\rho)\frac{\left(\frac12\right)_k}{z^k},\quad z\to\infty,\quad 
\tfrac12\le i\sigma_\pm<1,
\end{equation} 
where $(a)_k=\Gamma(a+k)/\Gamma(a)$ is the Pochhammer symbol. 
The first coefficients are
\begin{equation}\label{eq:asy03}
u_0(\rho)=-\frac{1}{\rho^2},\quad u_1(\rho)=\frac{\rho^2-2}{2\rho^4},\quad u_2(\rho)=-\frac{3\rho^4-4\rho^2+8}{8\rho^6}.
\end{equation}

\subsection{A uniform  asymptotic expansion}\label{sec:uniform}

Next, we want to obtain an expansion that is valid for small  positive values of  $ \vert\sigma_\pm\vert$. As explained in \cite[Chapter~21 and Part~7]{Temme:2015:AMI} and several  papers  \cite{Gil:2014:TAN}, \cite{Gil:2019:TCI} and \cite{Gil:2022:NAR}, we can use  the complementary error function to obtain uniform approximations. We have the integral representations (for properties of the error functions we refer to \cite{Temme:2010:ERF})
\begin{equation}\label{eq:asy04}
\begin{array}{r@{\,}c@{\,}l}
\erfc(z)&=&\dsp{\frac{2}{\sqrt{\pi}}\int_{z}^\infty e^{-t^2}\,dt,\quad z\in\CC,}\\[8pt]
w(z)&=&\dsp{\frac{1}{\pi i}\int_{-\infty}^{\infty} e^{-t^{2}}\frac{dt}{t-z}=
e^{-z^2}\erfc(-iz),\quad \Im z>0.}
\end{array}
\end{equation} 

In our analysis we need the function
\begin{equation}\label{eq:asy07}
V(\rho,z)=\int_{-\infty}^\infty e^{-z\sigma^2}\frac{d\sigma}{\sigma^2-\rho^2}=\frac{\pi }{i\rho}w\left(\rho\sqrt{z}\right)=
\frac{\pi }{i\sigma_\pm}e^{-z\sigma_\pm^2}
\erfc\left(i\sigma_\pm\sqrt{z}\right),
\end{equation} 
where we used $\rho=\sigma_\pm$ with $\sigma_\pm$ defined in \eqref{eq:prep07}. Both $\sigma_\pm$ satisfy $\Im \sigma_\pm<0$ (see also the third line of \eqref{eq:trans13}).

Next we consider $U(\sigma_\pm,z)$ defined in \eqref{eq:prep10}. We specify the role of the poles  by writing:
\begin{equation}\label{eq:asy11}\
\frac{1}{(\sigma^2-\rho^2)\sqrt{1+\sigma^2}}=g(\sigma,\rho)+\frac{1}{(\sigma^2-\rho^2)\sqrt{1+\rho^2}},\quad \rho=\sigma_\pm,
\end{equation} 
where
\begin{equation}\label{eq:asy12}
\begin{array}{r@{\,}c@{\,}l}
g(\sigma,\rho)&=&\dsp{\frac{1}{\sigma^2-\rho^2}\left(\frac{1}{\sqrt{1+\sigma^2}}-\frac{1}{\sqrt{1+\rho^2}}\right)}\\[8pt]
&=&\dsp{-\frac{1}{\sqrt{1+\sigma^2}\,\sqrt{1+\rho^2}\left(\sqrt{1+\sigma^2}+\sqrt{1+\rho^2}\right)}.}
\end{array}
\end{equation} 
This gives for the function  $U(\rho,z)$ defined in \eqref{eq:prep10}
\begin{equation}\label{eq:asy13}
\begin{array}{r@{\,}c@{\,}l}
U(\rho,z)&=&\dsp{\int_{-\infty}^\infty e^{-z\sigma^2}\left(\frac{1}{(\sigma^2-\rho^2)\sqrt{1+\rho^2}}+g(\sigma,\rho)\right)\,d\sigma}\\[8pt]
&=&\dsp{\frac{1}{\sqrt{1+\rho^2}}\,V(\rho,z)+\int_{-\infty}^\infty e^{-z\sigma^2} g(\sigma,\rho)\,d\sigma}\\[8pt]
&=&\dsp{-\frac{\pi i}{\rho\sqrt{1+\rho^2}}e^{-z\rho^2}\erfc\left(i\rho\sqrt{z}\right)+\int_{-\infty}^\infty e^{-z\sigma^2} g(\sigma,\rho)\,d\sigma}.
\end{array}
\end{equation} 

The function $g(\sigma,\rho)$ is analytic for $\vert\sigma\vert <1$ and we  expand
\begin{equation}\label{eq:asy14}
g(\sigma,\rho)=g(0,\rho)\sum_{k=0}^\infty c_k(\rho) \sigma^{2k},\quad g(0,\rho)=-\frac{1}{\sqrt{1+\rho^2}\left(1+\sqrt{1+\rho^2}\right)},
 \end{equation} 
 to obtain the asymptotic expansion
\begin{equation}\label{eq:asy15}
U(\rho,z)\sim-\frac{\pi i}{\rho\sqrt{1+\rho^2}}e^{-z\rho^2}\erfc\left(i\rho\sqrt{z}\right)+g(0,\rho)\sqrt{\frac{\pi}{ z}}\sum_{k=0}^\infty \frac{d_k(\rho)}{z^k},
\end{equation} 
where
\begin{equation}\label{eq:asy16}
d_k(\rho)= c_k(\rho) \left(\tfrac12\right)_k,\quad k=0,1,2,\cdots.
\end{equation} 

From the second line in \eqref{eq:asy12} we see that the computation of  the coefficients $c_k$ can be done by multiplying the Maclaurin series of $1/\sqrt{1+\sigma^2}$ and that of $1/\left(\sqrt{1+\sigma^2}+\sqrt{1+\rho^2}\right)$. For large values of $k$  both coefficients of these expansions are of small algebraic order of $k$, but  the coefficients  $d_k=c_k\left(\frac12\right)_k=c_k\Gamma\left(k+\frac12\right)/\Gamma\left(\frac12\right)$ are of factorial order.

After these preparations we have the asymptotic expansions for the functions $F^\pm(x;\alpha,\beta,\mu,\delta)$ that we defined in \eqref{eq:prep10}
\begin{equation}\label{eq:asy17}
F^\pm(x;\alpha,\beta,\mu,\delta)\sim \frac{e^{z\sigma_+^2}\sin(\nu\mp\tau)}{4\pi} 
\Bigl(
-\frac{\pi ie^{-z\sigma_\pm^2}}{\rho\sqrt{1+\rho^2}} \erfc\left(i\sigma_\pm\sqrt{z}\right)+g(0)\sqrt{\frac{\pi}{ z}}
\sum_{k=0}^\infty \frac{d_k(\sigma_\pm)}{z^k}
\Bigr).
\end{equation}

Let's  simplify a few expressions. We have  (see \eqref{eq:prep11} and\eqref{eq:asy14})
\begin{equation}\label{eq:asy18}
\begin{array}{r@{\,}c@{\,}l}
\dsp{\frac{\sin(\nu\mp\tau)}{\rho\sqrt{1+\rho^2}}}&=&\dsp{\frac{\sin(\nu\mp\tau)}{-i\sin\left(\frac12(\nu\mp\tau)\right)\cos\left(\frac12(\nu\mp\tau)\right)}=2i,}\\[8pt]
\sin(\nu\mp\tau)g(0,\sigma_\pm)
&=&\dsp{-\frac{\sin(\nu\mp\tau)}{\cos\left(\tfrac12(\nu\mp\tau)\right)\left(1+\cos\left(\tfrac12(\nu\mp\tau)\right)\right)}}\\[8pt]
&=&\dsp{-\frac{\sin\left(\frac12(\nu\mp\tau)\right)}{\cos^2\left(\frac14(\nu\mp\tau)\right)}=-2\tan\left(\tfrac14(\nu\mp\tau)\right),}\\[8pt]
z(\sigma_+^2-\sigma_-^2)&=&2\alpha\omega\sin(\nu)\sin(\tau)=2\gamma\delta.
\end{array}
\end{equation} 

In summary, we have the following theorem. The parameters $\nu$, $\tau$, $\sigma_\pm$ and $z$ 
are given in~\eqref{eq:prep11}.
\begin{theorem}\label{thm:thm02}
The functions $F^\pm(x;\alpha,\beta,\mu,\delta)$ composing the normal inverse Gaussian cumulative distribution $F(x;\alpha,\beta,\mu,\delta)$  as written  in \eqref{eq:prep10}  of Theorem~\ref{thm:thm01} have the asymptotic expansions
\begin{equation}\label{eq:asy19}
F^\pm(x;\alpha,\beta,\mu,\delta)\sim \frac{e^{z\sigma_+^2}}{4\pi} 
\Bigl(
2\pi e^{-z\sigma_\pm^2}\erfc\left(i\sigma_\pm\sqrt{z}\right)-
2\tan\left(\tfrac14(\nu\mp\tau)\right)\sqrt{\frac{\pi}{ z}}
\sum_{k=0}^\infty \frac{d_k(\sigma_\pm)}{z^k}
\Bigr),
\end{equation} 
or
\begin{equation}\label{eq:asy20}
\begin{array}{r@{\,}c@{\,}l}
F^+(x;\alpha,\beta,\mu,\delta)&\sim& \dsp{ \tfrac12\erfc(\zeta_+) - e^{-\zeta_+^2}
\frac{\tan\left(\frac14(\nu-\tau)\right)}{2\sqrt{\pi z}}
\sum_{k=0}^\infty \frac{d_k(\sigma_+)}{z^k}
,}\\[8pt]
F^-(x;\alpha,\beta,\mu,\delta)&\sim& \dsp{\tfrac12e^{2\gamma\delta}
\erfc(\zeta_-)-e^{-\zeta_+^2}\frac{\tan\left(\tfrac14(\nu+\tau)\right)}{2\sqrt{\pi z}}
\sum_{k=0}^\infty \frac{d_k(\sigma_-)}{z^k},}
\end{array}
\end{equation} 
where $2\gamma\delta=\zeta_-^2-\zeta_+^2$ and
\begin{equation}\label{eq:asy21}
\begin{array}{r@{\,}c@{\,}l}
\zeta_+&=&i\sigma_+\sqrt{z}=\sin\left(\frac12(\nu-\tau)\right)\sqrt{z}=\sign(x_0-x)\sqrt{\alpha\omega-\beta\xi-\gamma\delta},\\[8pt]
\zeta_-&=&i\sigma_-\sqrt{z}=\sin\left(\frac12(\nu+\tau)\right)\sqrt{z}=\sqrt{\alpha\omega-\beta\xi+\gamma\delta}.
\end{array}
\end{equation}
The expansions are valid for large values of $z$ and $\vert\sigma_\pm\vert\le \frac12$.
\end{theorem}

\vspace{3mm}

For the complementary function $G(x;\alpha,\beta,\mu,\delta)$ defined in \eqref{eq:intro02} we have
\begin{corollary}\label{cor:cor01}
The asymptotic expansion of the function $G(x;\alpha,\beta,\mu,\delta)$ follows from
\begin{equation}\label{eq:asy22}
G(x;\alpha,\beta,\mu,\delta)=G^+(x;\alpha,\beta,\mu,\delta)-F^-(x;\alpha,\beta,\mu,\delta)
\end{equation}
and 
\begin{equation}\label{eq:asy23}
G^+(x;\alpha,\beta,\mu,\delta)\sim \tfrac12\erfc(-\zeta_+) + e^{z\sigma_+^2}
\frac{\tan\left(\frac14(\nu-\tau)\right)}{2\sqrt{\pi z}}
\sum_{k=0}^\infty \frac{d_k(\sigma_+)}{z^k}.
\end{equation} 
\end{corollary}

\begin{figure}[tb]
\vspace*{0.0cm}
\begin{center}
\begin{minipage}{5.0cm}
        \includegraphics[width=5.0cm]{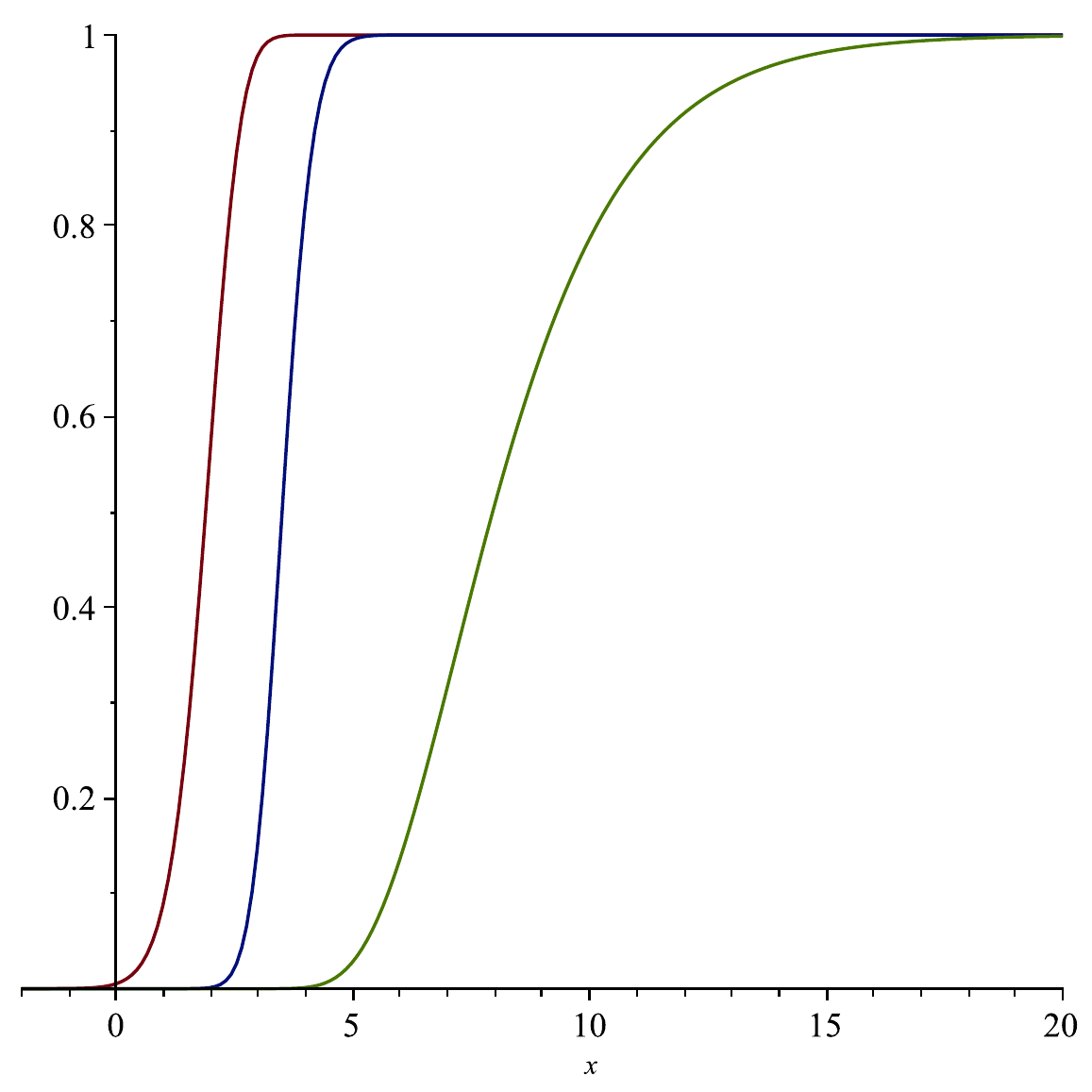} 
\end{minipage}
\hspace*{2cm}
\begin{minipage}{5.5cm}
        \includegraphics[width=5.0cm]{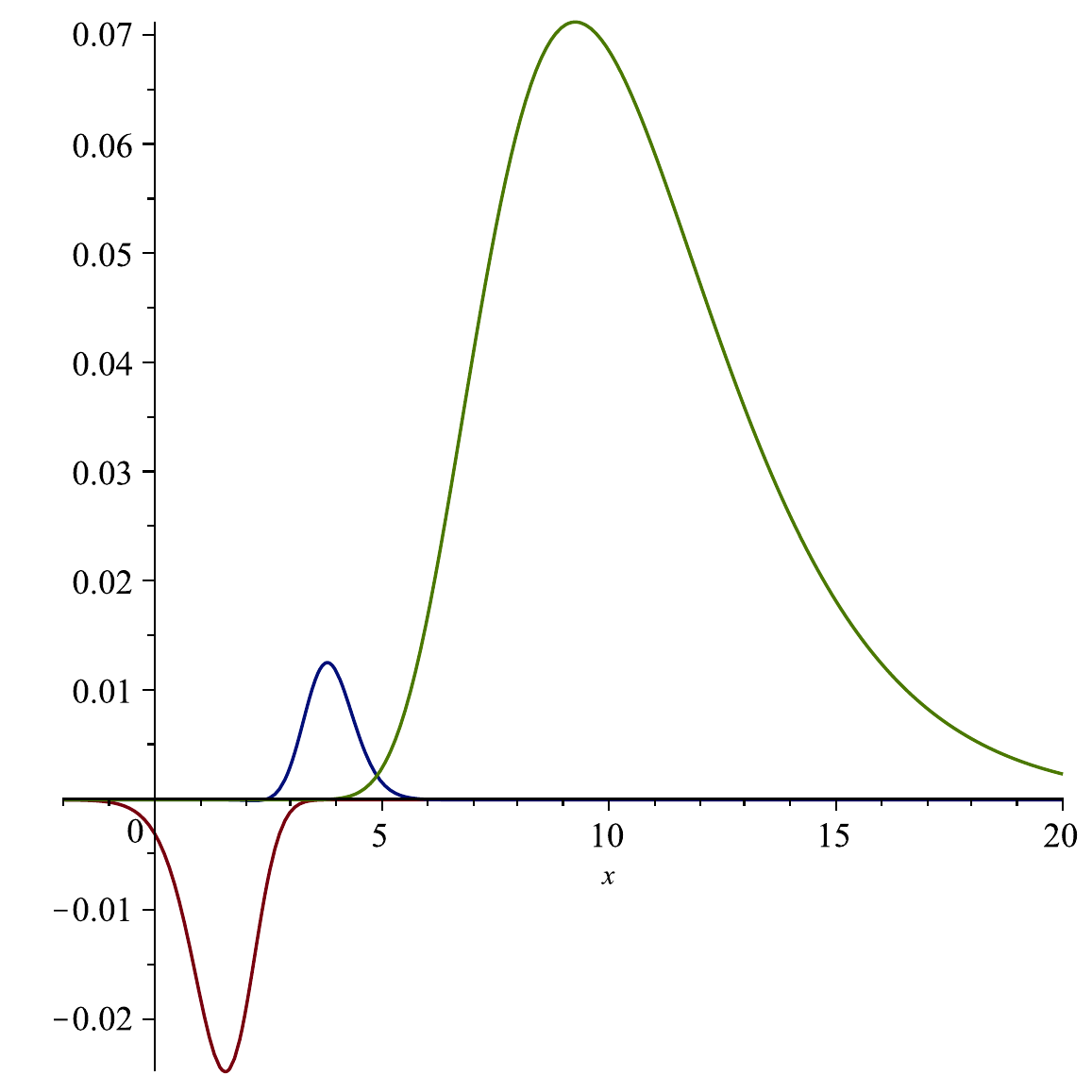} 
  \end{minipage}
\end{center}
\caption{Graphs obtained using the asymptotic approximations,  for $\alpha=8$, $\mu=3$, $\delta=2$, $0 \le x \le 20$, and three values of $\beta$: left $\beta=-4$, middle $\beta=2$, right $\beta=7.5$. In the left figure we see graphs of the function $F(x;\alpha,\beta,\mu,\delta)$, in the right one we see $F^-(x;\alpha,\beta,\mu,\delta)$ for the same parameters.  }
\label{fig:fig01}
\end{figure}


We make a few observations.
\begin{enumerate}
\item
The square root forms of $\zeta_\pm$ in \eqref{eq:asy21} follow from the relations in \eqref{eq:prep11}. The expression $\alpha\omega-\beta\xi-\gamma\delta$ for $\zeta_+$ is a  convex functions of $\xi$ with a double zero at $\xi_0$, the transition point. This follow from calculating the first terms in the Taylor expansion of $\zeta_+$ as function of $\xi$ at  $\xi_0$.
\item
The argument $\zeta_+$ of $\erfc(\zeta_+)$  in the first line of \eqref{eq:asy20} vanishes when
$\nu=\tau$, that is, when $x=x_0$; see the earlier discussion at equation~\eqref{eq:trans12}. 
In the first paragraph of Section~\ref{sec:prep} we assumed that $\nu>\tau$ (or $x <x_0$), since we have started with \eqref{eq:trans10}. This corresponds to $\zeta_+>0$ in $\erfc(\zeta_+)$.
\item
However, $\erfc(\zeta_+)$ allows us to take $\nu <\tau$ (or $x>x_0$), and the  change from  $\nu>\tau$ to $\nu<\tau$ in $\erfc(\zeta_+)$ is smooth and analytic. Also other elements in the series expansion of the function  $F^+(x;\alpha,\beta,\mu,\delta)$ remain well-defined, and $F^+(x;\alpha,\beta,\mu,\delta)$   changes from values smaller than $\frac12$ ($x\le x_0$ or $\nu>\tau$) to values larger than $\frac12$ ($x\ge x_0$ or $\nu < \tau$).

\item
When $\nu < \tau$ w can repeat the analysis with the representation of the complementary function $G(x;\alpha,\beta,\mu,\delta)$, starting with the first line of \eqref{eq:trans09}. The result will be the same as in Corollary~\ref{cor:cor01}.

\item
 A notable point is that the  asymptotic approximations can be used for negative values of $\xi=x-\mu$, while the starting representation in \eqref{eq:trans01} is only valid for $\xi\ge0$. 
 See for $\xi<0$ also the representation in \eqref{eq:intro06}.

\end{enumerate}

\vspace{3mm}

The first coefficients $d_k(\sigma_\pm)$ are given by
\begin{equation}\label{eq:asy24}
\begin{array}{r@{\,}c@{\,}l}
  d_0(\sigma_\pm)&=&\dsp{ 1,\quad
  d_1(\sigma_\pm)=  -\frac{w+2}{4 (w+1)},\quad
  d_2(\sigma_\pm)=  \frac{3 (3 w^2+9 w+8)}{32(w+1)^2},}\\[8pt]
  d_3(\sigma_\pm)&=&\dsp{  -\frac{ 15(5w^3+20w^2+29w+16)}{128(w+1)^3},}\\[8pt]
  d_4(\sigma_\pm)&=&\dsp{\frac{105(35w^4+175w^3+345w^2+325w+128)}{2048(w+1)^4}.}
\end{array}
\end{equation} 
Here,  $w=\sqrt{1+\sigma_\pm^2}$. Hence, $w=\cos\left(\frac12(\nu\mp\tau)\right)$. 
In the Appendix we describe a recursive method for the symbolic evaluation of these coefficients, together with a short Maple code.

In Figure~\ref{fig:fig01} we show graphs obtained using the asymptotic approximations  in \eqref{eq:asy20},  for $\alpha=8$, $\mu=3$, $\delta=2$, $0 \le x \le 20$, and three values of $\beta$: left $\beta=-4$, middle $\beta=2$, right $\beta=7.5$. In the left figure we see 3 graphs of the function $F(x;\alpha,\beta,\mu,\delta)$, in the right one we see 3 graphs of  $F^-(x;\alpha,\beta,\mu,\delta)$ for the same parameters. We see that the main contributions are coming from $F^+(x;\alpha,\beta,\mu,\delta)$, but for numerical computations  $F^-(x;\alpha,\beta,\mu,\delta)$ cannot be neglected.

The corresponding transition values $x_0$, function values of $F(x_0;8,\beta,3,2)$, and values of the large parameter $z=2\alpha\omega$ are given in Table~\ref{tab:table01}.
In addition, we show the absolute errors in the computation of the asymptotic expansions relative to more accurate computations of the $F$ function. The computations of the asymptotic approximations are performed with Maple, $Digits=8$, with terms up to and including $k=5$ in the asymptotic expansions in \eqref{eq:asy20}. The absolute errors shown illustrate the quality of the asymptotic approximation with only 6 terms. In a future paper, we plan to perform more extensive tests to assess the accuracy of the asymptotic expansions for different regions of the parameters.

\begin{table}[th]\renewcommand{\arraystretch}{1.19}\label{tab:tab01}
$$
\begin{array}{rllllll}
\hline
 \beta & \quad\quad x_0 & \quad\quad\quad F &\quad\quad z& {\rm absolute \ errors} \\
 \hline
 -4.0  &    1.845299462 & 0.473833601   & 36.95041722 &3.1\times 10^{-08}\\
 2.0   &    3.516397780 & 0.512385772 & 33.04945788 & 1.7\times 10^{-09} \\
 7.5   &     8.388159062  & 0.575900502 & 91.95791466 &4.7\times 10^{-10}&
\end{array}
$$

\vspace{-2mm}

\caption{For the values of $\beta$ used in Figure~\ref{fig:fig01} we give the transition values $x_0$, the function values $F(x_0;\alpha,\beta,\mu,\delta)$ for $\alpha=8$, $\mu=3$, $\delta=2$,  the large parameter $z=2\alpha\omega$, and absolute errors in the shown values of the $F$-function for these parameters. }
\label{tab:table01}
\end{table}

\section{Appendix}\label{sec:app}
We can obtain the coefficients $c_k$ of the expansion \eqref{eq:asy14} of the function $g(\sigma,\rho)$ given in \eqref{eq:asy12} using Maple's  procedure {\it taylor}. To describe an algorithm without this procedure, and without Maple's procedure {\it pochhammer} we use a recursive method for the symbolic evaluation of the coefficients $c_k$. 
Finally, we provide a Maple procedure {\it dkproc} for this method. Similar code can be written for the evaluation of the coefficients $u_k(\rho)$ of the Maclaurin series in \eqref{eq:asy01}.

We use the second line of equation \eqref{eq:asy12} and write the expansion in \eqref{eq:asy14} in the form
\begin{equation}\label{eq:app01}
-g(0,\rho)\left(w(1+\sigma^2)+w^2\sqrt{1+\sigma^2}\right)\sum_{k=0}^\infty  c_k(\rho)\sigma^{2k}=1,\quad w=\sqrt{1+\rho^2},
\end{equation} 
or
\begin{equation}\label{eq:app02}
-g(0,\rho)\sum_{k=0}^\infty  b_k(\rho)\sigma^{2k}\sum_{k=0}^\infty  c_k(\rho)\sigma^{2k}=1,
\end{equation} 
where
\begin{equation}\label{eq:app03}
c_0=1,\quad b_0=w+w^2,\quad b_1=w+\tfrac12 w^2, \quad b_k=w^2(-1)^k\frac{(-\frac12)_k}{k!},\quad k\ge 2.
\end{equation} 
The next step is multiplying the two series:
\begin{equation}\label{eq:app04}
-g(0,\rho)\sum_{k=0}^\infty  a_k(\rho)\sigma^{2k}=1,\quad a_k=\sum_{j=0}^k b_j (\rho)c_{k-j}(\rho).
\end{equation} 
All coefficients $a_k(\rho)$ should vanish, except $a_0(\rho)=b_0(\rho)c_0(\rho)$, and we have using \eqref{eq:asy14} the equation $-g(0,\rho)b_0(\rho)c_0(\rho)=1$, which yields $c_0(\rho)=1$. For the other $c_k(\rho)$ we obtain the recursive relation
\begin{equation}\label{eq:app05}
c_k(\rho)=-\
\frac{1}{b_0(\rho)}\sum_{j=1}^{k} b_j(\rho) c_{k-j}(\rho),\quad k=1,2,3,\cdots.
\end{equation} 
With these coefficients we can obtain $d_k(\rho)=\left(\frac12\right)_k c_k(\rho)$.

In the following Maple code we describe an algorithm to compute $d_k(\rho)$.
\begin{verbatim}
restart;
dkproc:= proc(kmax, dk) local bk, ck, w, j, k, s, p;
# To compute d_k, see (4.24); dk is an output parameter.
  bk[0]:= w+w^2; bk[1]:= w+w^2/2; bk[2]:= -w^2/8; 
  ck[0]:= 1; dk[0]:= 1; 
  ck[1]:= normal(-bk[1]/bk[0]); 
  p:= 1/2; dk[1]:= p*ck[1];
  # p = (1/2)_1; a Pochhammer value to start a recursion.
  for k from 2 to kmax do 
    # bk[k] is a binomial coefficient, generated by recursion.
    bk[k+1]:= -(k-1/2)*bk[k]/(k+1); 
    s:= 0; 
    for j from 1 to k do s:= normal(s + bk[j]*ck[k-j]) od;  
    ck[k]:= normal(-s/bk[0]);  p:= p*(k-1/2); 
    #This makes p = (1/2)_k;
     dk[k]:= p*ck[k];
  od;
  kmax
end;
# For example:
kmax:= 5;
dkproc(kmax, dk);
for k from 0 to kmax do print(k,dk[k]) od;
\end{verbatim}

\section*{Acknowledgements}
The author thanks the reviewers for their constructive comments and helpful suggestions. One reviewer verified  the paper in detail, noted in particular that the transition value $x_0$ coincides with the mean of the normal Gaussian distribution (see below (\eqref{eq:trans12})), noted the typo in the first line of (3.10), and provided absolute errors for Table~\ref{tab:tab01}. The other reviewer suggested, among other things, including a Maple code for the coefficients $d_k(\rho)$ of the asymptotic expansions (see the new Appendix), and provided a 31-page transcript of a ChatGPT session showing the verification of all steps in the paper, including transformations of variables, shifts of integration paths into the complex plane across the poles, and so on.\\
The author thanks Amparo Gil, Guillermo Navas-Palencia and Javier Segura for their advices at the start of this project.\\
The author acknowledges financial support from project PID2021-
127252NB-I00 funded by MCIN/AEI/10.13039/501100011033/ FEDER, UE.\\

\bibliographystyle{plain}
\bibliography{Nig}

\end{document}